\documentclass[12pt,leqno,draft]{article}

\newtheorem{theorem}{Theorem}
\newtheorem{lemma}[theorem]{Lemma}
\newtheorem{proposition}[theorem]{Proposition}
\newtheorem{definition}[theorem]{Definition}
\newtheorem{corollary}[theorem]{Corollary}

\newcommand{\begintheorem}{\addtocounter{equation}{1}\begin{theorem}}
\newcommand{\beginlemma}{\addtocounter{equation}{1}\begin{lemma}}
\newcommand{\beginproposition}{\addtocounter{equation}{1}\begin{proposition}}
\newcommand{\begindefinition}{\addtocounter{equation}{1}\begin{definition}}
\newcommand{\begincorollary}{\addtocounter{equation}{1}\begin{corollary}}

\begin{document}

\title{Elements of harmonic analysis, 5}

\author{Stephen William Semmes	\\
	Rice University		\\
	Houston, Texas}

\date{}

\maketitle

	Fix a positive integer $n$, and let $\mathcal{S}({\bf R}^n)$
denote the vector space of complex-valued smooth functions $f(x)$
on ${\bf R}^n$ such that $f$ and all of its derivatives are bounded
functions on ${\bf R}^n$, and moreover they are rapidly decreasing
in the sense that
\begin{equation}
\label{x^alpha frac{partial^{d(beta)}}{partial x^beta} f(x)}
	x^\alpha \frac{\partial^{d(\beta)}}{\partial x^\beta} f(x)
\end{equation}
are also bounded functions on ${\bf R}^n$ for all multi-indices
$\alpha$, $\beta$, where $d(\beta)$ denotes the degree of $\beta$,
which is to say the sum of its components.  Thus $\mathcal{S}({\bf
R}^n)$ contains the smooth functions on ${\bf R}^n$ with compact
support as a linear subspace.  Also, $\mathcal{S}({\bf R}^n)$ contains
the Gaussian functions $\exp (- A(x) \cdot x)$, where $A$ is a
positive-definite symmetric linear transformation on ${\bf R}^n$, and
$x \cdot y$ is the usual inner product defined by
\begin{equation}
	x \cdot y = \sum_{j=1}^n x_j \, y_j.
\end{equation}

	The space $\mathcal{S}({\bf R}^n)$ is commonly known as the
\emph{Schwartz class} of rapidly decreasing smooth functions on ${\bf
R}^n$.  If $f_1$, $f_2$ are functions in the Schwartz class, then the
sum $f_1 + f_2$ and the product $f_1 \, f_2$ also lie in the Schwartz
class.  If $f \in \mathcal{S}({\bf R}^n)$, then each translate $f(x -
a)$ lies in $\mathcal{S}({\bf R}^n)$ too, and derivatives of $f$ of
all orders are also elements of $\mathcal{S}({\bf R}^n)$.
Furthermore, the product of $f$ with a polynomial, or a function of
the form $\exp (i \, b \cdot x)$, $b \in {\bf R}^n$, lie in
$\mathcal{S}({\bf R}^n)$ too.

	Let $f(x)$ be a function in the Schwartz class on ${\bf R}^n$.
Because $f$ is rapidly decreasing, it follows that $f$ is integrable
on ${\bf R}^n$.  Thus we can define the Fourier transform of $f$
as usual by
\begin{equation}
	\widehat{f}(\xi) 
	   = \int_{{\bf R}^n} f(x) \, \exp (- 2 \pi i \, \xi \cdot x) \, dx.
\end{equation}
Just as for any integrable function on ${\bf R}^n$, or any finite
measure, the Fourier transform is a bounded continuous function on
${\bf R}^n$, and even uniformly continuous.  For functions in the
Schwartz class we have a much stronger result, namely that the Fourier
transform also lies in the Schwartz class.

	Indeed, $\widehat{f}$ is a bounded continuous function on
${\bf R}^n$, and for each multi-index $\alpha$ we have that $x^\alpha
\, f(x)$ is an integrable function on ${\bf R}^n$, and hence its
Fourier transform is also a bounded continuous function on ${\bf
R}^n$.  Using this one can show that $\widehat{f}$ is smooth, and that
derivatives of all orders of $\widehat{f}$ are bounded.  For each
multi-index $\beta$ we have that $\xi^\beta \, \widehat{f}(\xi)$ is a
constant times the Fourier transform of $(\partial^{d(\beta)} /
\partial x^\beta) f(x)$, which implies that $\xi^\beta \,
\widehat{f}(\xi)$ is the Fourier transform of an integrable function,
and is therefore bounded.  Similarly one can show that derivatives of
all orders of $\widehat{f}$ decay rapidly, so that $\widehat{f} \in
\mathcal{S}({\bf R}^n)$.

	If $\phi$ is an integrable function on ${\bf R}^n$, then
the inverse Fourier transform of $\phi$ is defined by
\begin{equation}
	{\check \phi}(y) 
  = \int_{{\bf R}^n} \phi(\eta) \, \exp (2 \pi i \, y \cdot \eta) \, d\eta.
\end{equation}
As for the Fourier transform, the inverse Fourier transform of an
integrable function, or more generally a finite measure, is a bounded
uniformly continuous function on ${\bf R}^n$.  If $\phi \in
\mathcal{S}({\bf R}^n)$, then ${\check \phi} \in \mathcal{S}({\bf
R}^n)$ too, for essentially the same reasons as for the Fourier
transform.  In general, if $f$ is a continuous integrable function on
${\bf R}^n$ such that $\widehat{f}$ is also integrable, then the
inverse Fourier transform of the Fourier transform of $f$ is equal to
$f$ again.  In particular, this applies to the case where $f$ is in
the Schwartz class, since the Fourier transform of $f$ is then in the
Schwartz class as well.

	Similarly, if $\phi$ is a continuous integrable function on
${\bf R}^n$ such that ${\check \phi}$ is integrable, then the Fourier
transform of ${\check \phi}$ is defined.  Under these conditions the
Fourier transform of the inverse Fourier transform of $\phi$ is equal
to $\phi$, and in particular this holds when $\phi \in
\mathcal{S}({\bf R}^n)$.

	To summarize a bit, the Fourier transform and inverse Fourier
transforms define linear mappings from the Schwartz class into itself
which are indeed inverses of each other.  It follows that the Fourier
transform and inverse Fourier transform are one-to-one linear mappings
of the Schwartz class onto itself.  Thus every element of the Schwartz
class arises as the Fourier transform of an element of the Schwartz
class.

	If $f_1$, $f_2$ are two elements of the Schwartz class on
${\bf R}^n$, then the convolution of $f_1$ and $f_2$ is defined by
\begin{equation}
	(f_1 * f_2)(x) = \int_{{\bf R}^n} f_1(y) \, f_2(x - y) \, dy.
\end{equation}
The Fourier transform of $f_1 * f_2$ is equal to the product
of the Fourier transforms of $f_1$, $f_2$.  Because $f_1$, $f_2$
are elements of the Schwartz class, so are their Fourier transforms.
The product of the Fourier transforms of $f_1$, $f_2$ therefore
lies in the Schwartz class, and one may conclude that $f_1 * f_2$
is an element of the Schwartz class too.

	For each function $f$ in the Schwartz class on ${\bf R}^n$ and
each pair $\alpha$, $\beta$ of multi-indices, define $\|f\|_{\alpha,
\beta}$ to be the supremum of the absolute value of (\ref{x^alpha
frac{partial^{d(beta)}}{partial x^beta} f(x)}) over all $x \in {\bf
R}^n$.  This defines a seminorm on $f$, which means that if $f_1$,
$f_2$ are elements of the Schwartz class, then $\|f_1 + f_2\|_{\alpha,
\beta}$ is less than or equal to the sum of $\|f_1\|_{\alpha, \beta}$
and $\|f_2\|_{\alpha, \beta}$, and that if $f$ is an element of the
Schwartz class and $c$ is a complex number, then $\|c \, f\|_{\alpha,
\beta}$ is equal to the product of the absolute value of $c$ and
$\|f\|_{\alpha, \beta}$.  Using this family of seminorms on the
Schwartz class we get a topology, which is generated by finite
intersections of balls defined with respect to the seminorms.

	It is perhaps more convenient to think of this topology in
terms of sequences.  Namely, if $\{f_j\}_{j=1}^\infty$ is a sequence
of functions in the Schwartz class and $f$ is another function in the
Schwartz class, then $\{f_j\}_{j=1}^\infty$ converges to $f$ in the
topology determined by the seminorms if and only if $\|f_j -
f\|_{\alpha, \beta} \to 0$ as $j \to \infty$ for all multi-indices
$\alpha$, $\beta$.  We can say that a sequence $\{f_j\}_{j=1}^\infty$
of functions in the Schwartz class is a Cauchy sequence if $\|f_j -
f_l\|_{\alpha, \beta} \to 0$ as $j, l \to \infty$ for all
multi-indices $\alpha$, $\beta$.  One can show that the Schwartz class
is complete in the sense that any Cauchy sequence in the Schwartz
class converges in the sense just defined to an element of the
Schwartz class.  The Fourier transform and inverse Fourier transform
are continuous linear transformations from the Schwartz class to
itself.

	The vector space operations of addition and scalar
multiplication define continuous mappings from the Schwartz class to
itself, as does multiplication of two elements of the Schwartz class.
Linear transformations on the Schwartz class defined by
differentiating a function in the Schwartz class to some order, or
multiplying it by some fixed polynomial, are also continuous.  The
Schwartz class on ${\bf R}^n$ is a basic and natural example of a
topological vector space which is a Fr\'echet space, which means that
the topology can be defined by a countable family of seminorms and
that the space is complete, while the topology is not determined by
any single norm.

	By a tempered distribution on ${\bf R}^n$ we mean a continuous
linear mapping from the Schwartz class on ${\bf R}^n$ into the complex
numbers.  Explicitly, a linear functional $\lambda$ on
$\mathcal{S}({\bf R}^n)$ is continuous if there is a positive real
number $C$ and a finite collection of pairs of multi-indices
$\alpha_1, \beta_1, \ldots, \alpha_l, \beta_l$ such that
$|\lambda(f)|$ is less than or equal to $C$ times the sum of
$\|f\|_{\alpha_j, \beta_j}$ over $1 \le j \le l$ for all functions $f$
in the Schwartz class.  A basic class of examples consists of linear
functionals $\lambda$ defined by setting $\lambda(f)$ equal to the
integral over ${\bf R}^n$ of the product of $f$ and a continuous
function $h$ of polynomial growth.  If $\lambda$ is a tempered
distribution and $\alpha$ is a multi-index one can define the
$\partial^{d(\alpha)} / \partial x^\alpha$ derivative of $\lambda$ as
a tempered distribution by saying that this distribution applied to a
function $f$ in the Schwartz class is the same as $(-1)^{d(\alpha)}$
times $\lambda$ applied to $(\partial^{d(\alpha)} / \partial x^\alpha)
f$, which is the same as what one would get by integration by parts if
$\lambda$ is defined by integration of $f$ times a smooth density $h$.
Similarly, one can define the Fourier transform $\widehat{\lambda}$ of
a tempered distribution by saying that $\widehat{\lambda}(f)$ is equal
to $\lambda(\widehat{f})$ for all $f$ in the Schwartz class.

\end{document}